\documentclass[10.5pt]{article} 
\usepackage{fouriernc}
\usepackage{extarrows}
\usepackage{array}

\usepackage{amsfonts,amsmath,amssymb,amsthm,amsxtra,graphicx,latexsym}
\numberwithin{equation}{section}

\def\sct{\section}
\def\ssct{\subsection}

\def\ba{\sla\begin{array}}\def\ea{\end{array}}     
\def\be{\begin{equation}}\def\ee{\end{equation}}   
\def\bgd{\begin{aligned}}\def\egd{\end{aligned}}   
    \def\lgd{\lt\{\bgd}\def\rgd{\egd\rt.}          
\def\bgt{\begin{gathered}}\def\egt{\end{gathered}} 

\def\bcc{\begin{center}}\def\ecc{\end{center}}
\def\bmc{\begin{multicols}}\def\emc{\end{multicols}}
\def\bmpg{\begin{minipage}}\def\empg{\end{minipage}}
\def\btb{\begin{tabular}}\def\etb{\end{tabular}}
\def\bfl{\begin{flushleft}}
\def\efl{\end{flushleft}}

\def\lb{\label}
\def\rf{\ref}

\def\<{\langle}
\def\>{\rangle}

\def\bs{\setminus}
\def\cd{\cdots}

\def\ept{\emptyset}
\def\ev{\equiv}

\def\fa{\forall}
\def\gsl{\geqslant}
\def\hra{{\hookrightarrow}}
\def\ift{\infty}

\def\lsl{\leqslant}

\def\lt{\left}\def\rt{\right}
\def\mrg{\mathring}
\def\mt{\mapsto}

\def\nev{\not\equiv}

\def\op{\oplus}

\def\pt{\partial}
\def\ra{{\rightarrow}}

\def\rap{\rightharpoonup}

\def\sbs{\subset}
\def\sps{\supset}

\def\ts{\times}
\def\udb{\underbrace}

\def\wtd{\widetilde}

\def\xra{\xrightarrow}

\def\aa{\alpha}
\def\bb{\beta}
\def\ga{\gamma}\def\Ga{\Gamma}
\def\dl{\delta}

\def\tht{\theta}

\def\lm{\lambda}

\def\sg{\sigma}

\def\om{\omega}

\def\R{\mathbb R}
\def\Z{\mathbb Z}

\def\cE{\mathcal E}

\def\cL{\mathcal L}

\def\cP{\mathcal P}
\def\rd{\textrm d}

\makeatletter
\def\ExtendSymbol#1#2#3#4#5{\ext@arrow 0099{\arrowfill@#1#2#3}{#4}{#5}}
\def\RightExtendSymbol#1#2#3#4#5{\ext@arrow 0359{\arrowfill@#1#2#3}{#4}{#5}}
\def\LeftExtendSymbol#1#2#3#4#5{\ext@arrow 6095{\arrowfill@#1#2#3}{#4}{#5}}
\makeatother

\newcommand\mra[2][]{\ExtendSymbol{~}{-}{\rightarrow}{#1}{#2}}

\DeclareSymbolFont{EUEX}{U}{euex}{m}{n}
\DeclareSymbolFont{euexlargesymbols}{U}{euex}{m}{n}
\DeclareMathSymbol{\intop}{\mathop}{euexlargesymbols}{"52}
\def\int{\intop\nolimits}

\DeclareSymbolFont{mathdesign}{OMX}{mdbch}{m}{n}
\DeclareMathSymbol{\intmd}{\mathop}{mathdesign}{"52}
\DeclareMathSymbol{\ointmd}{\mathop}{mathdesign}{"49}
\def\int{\intop\nolimits}

\DeclareSymbolFont{ugmL}{OMX}{mdugm}{m}{n}
\SetSymbolFont{ugmL}{bold}{OMX}{mdugm}{b}{n}
\DeclareMathAccent{\wideparen}{\mathord}{ugmL}{"F3}

\def\QEDclosed{\mbox{\rule[0pt]{0.7ex}{1.3ex}}} 

\def\QED{\QEDclosed} 

\def\dsty{\displaystyle}
\def\hs{\hspace}    

\def\nid{\noindent} 
\def\nn{\nonumber}  
\def\sla{\setlength{\arraycolsep}{2pt}} 
\def\sss{\scriptscriptstyle}

\def\crit{\text{Crit}}
\def\ncrit{^\#\crit}

\def\ncP{^\#\cP}

\newtheoremstyle{plain}{0.6ex}{0.6ex}{}{}{}{}{0.5em}{\bf\thmname{#1}~\thmnumber{#2}.~\thmnote{(#3)}}
\newtheoremstyle{definition}{0.6ex}{0.6ex}{}{}{}{}{0.5em}{\bf\thmname{#1}~\thmnumber{#2}.~\thmnote{(#3)}}

\theoremstyle{plain}
\newtheorem{thm}{Theorem}
\newtheorem{prop}{Proposition}[section]
\newtheorem{cor}{Corollary}

\theoremstyle{definition}
\newtheorem{Def}{Definition}
\newtheorem{rem}{Remark}

\def\bT{\begin{thm}}\def\eT{\end{thm}}
\def\bP{\begin{prop}}\def\eP{\end{prop}}
\def\bcr{\begin{cor}}\def\ecr{\end{cor}}
\def\bD{\begin{Def}}\def\eD{\end{Def}}
\def\bR{\begin{rem}}\def\eR{\end{rem}}
\def\pf{\noindent{\it Proof.~~}}

\def\rT{Theorem~\rf}
\def\rP{Proposition~\rf}

\allowdisplaybreaks

\begin{document}
\title{
\vspace{0.5in} {\bf\Large Multiple Forced rotations\\ for the N-pendulum Equation}}
\author{{\bf\large Hui Qiao\vspace{1mm}}\\
{\it\small School of Mathematics and Statistics}\\
{\it\small Wuhai University, Wuhan 430072},\\
{\it\small The People's Republic of China}\\
{\it\small e-mail: qiaohuimath@whu.edu.cn}}\vspace{2mm}
\date{
}
\maketitle
\begin{center}
{\bf\small Abstract}

\vspace{3mm} \hspace{.05in}\parbox{4.5in}
{\small In this paper, we consider the planar forced $N$-pendulum equation.
Multiple rotational solutions are obtained.}
\end{center}

\section{Introduction and main results}

Let $m_i,\ell_i$ be positive constants, $i=1,2,\cd,N$, and
\be\aa_j=\sum^N_{s=j}m_s,~\bb_j=\aa_j\ell_j,~j=1,2,\cd,N.\lb{eq:defabj}\ee
The motion of the coplanar $N$-pendulum is governed by the Lagrangian
\be L(q,p)=\frac{1}{2}A(q)p\cdot p+V(q),\lb{eq:defL}\ee
where
\[q=(q_1,q_2,\cd,q_{\sss N}),~~p=(p_1,p_2,\cd,p_{\sss N})\in\R^N,\]
\be A(q)p\cdot p=\sum_{1\lsl i\lsl N}\aa_i\ell^2_ip^2_i
+\sum_{1\lsl i<j\lsl N}2\aa_j\ell_i\ell_j\cos(q_i-q_j)p_ip_j,\lb{eq:defA(q)pp}\ee
\be V(q)=g\sum^N_{j=1}\bb_j\cos(q_j),\lb{eq:defV}\ee
In the following, we add an periodic forcing term $f(t)=(f_1(t),\cd,f_{\sss N}(t))$.
The corresponding system of equations are
\be \frac{\rd}{\rd t}L_p(q,\dot q)-L_q(q,\dot q)=f(t),\lb{eq:Lsys}\ee
Periodic motions for the N-pendulum are solutions of \eqref{eq:Lsys} satisfying
\be q(t+T)-q(t)=2\pi v,~v\in\Z^N,~\fa t\in\R.\lb{eq:pslt}\ee

In \eqref{eq:pslt}, if $v=0$, solutions are oscillatory. Existences and
multiplicities of oscillations for planar double, triple and $N$-pendulum have
been studied in
\cite{CFS87,FW89,CLZ90,Ta90,RP01} and references therein. In \cite{pR88},
Rabinowitz showed the existence of at least $N+1$ $2\pi$-periodic forced
oscillations for more general Lagrangian functions.

In this paper, we consider solutions with
\[v=(v_1,v_2,\cd,v_{\sss N})\ne0,\]
which are called {\it rotational solutions} and non-contractible on $T^N$. Note that
if $q(t)$ is a solution, then $q(kT+t)-q(t)\ev kv$. We restrict choices of rotational vectors to
$\Z^N_1$ defined as follows.

\bD A rotational vector $v=(v_1,v_2,\cd,v_{\sss N})\in\Z^N\bs\{0\}$ is called {\it prime}, if
one of its coordinates equals to $1$ while others are zero, or two of them are relatively
prime. Denotes by $\Z^N_1$ the set of all prime rotational vectors in $\Z^N$.\eD

For $v\in\Z^N_1$, let
\be q(t)=x(t)+\frac{2\pi vt}{T},~\dot q(t)=\dot x(t)+\frac{2\pi v}{T},\lb{eq:q=x+tv}\ee
and define a functional $\cL:E=W^{1,2}(S^1,\R^N)\ra\R$ by
\be \cL(x)=\int^T_0\big(L(q(t),\dot q(t))+f(t)\cdot q(t)\big)\rd t.\lb{eq:defcA}\ee
Its critical points give rotational solutions. When we study multiplicities, the following symmetries are involved:

{\it $\Z^N$-symmetry}. Note that $E$ can be orthogonally decomposed as
\be E=E^0\op\wtd E,\lb{eq:E=E0+E'}\ee
where $E^0$ denotes the subspace of constant functions and $\wtd E$ denotes the subspace
of functions with zero mean value. Since the function $L$ is $2\pi$-periodic in $q_i$
for all $1\lsl i\lsl N$ and $f\in\wtd E$, the functional $\cL$ defined by \eqref{eq:defcA}
is also $2\pi$-periodic in all $q_i$s. Hence it can be defined on
\be\cE=E/(2\pi\Z^N)=T^N\ts\wtd E.\lb{eq:cE=TNxE'}\ee

{\it $S^1$-symmetry}. If $f\ev0$, there is an $S^1$-action
involved in the problem. Given $v\in\Z^N_1$ and $T>0$, let
\be\tht\cdot x(t)=x(t+\tht)+\frac{2\pi v}{T}\tht,\lb{eq:L.s1act}\ee
where we use the representation of $S^1$ by $\dsty\lt\{\exp\lt(\frac{2\pi\sqrt{-1}}{T}\tht\rt)~\bigg|~\tht\in[0,T]\rt\}$.
Then $\cL$ is $S^1$-invariant. Moreover, this action is free.

Denote by $\cP(v,T)$ the set of distinct rotational solutions of \eqref{eq:Lsys} with
given rotational vector $v$ and period $T$.

\bT\lb{thm:Lsys.rot}\it Assume that $f$ satisfies
\be f\in L^2(\R/(T\Z),\R^N),~\int^T_0f(t)\rd t=0.\lb{eq:intf=0}\ee
For every $v\in\Z^N_1$ and $T>0$, we have\\
(i) If $f\nev0$, then
\[\ncP(v,T)\gsl N+1.\]
If all solutions are nondegenerate, then $\ncP(v,T)\gsl 2^N$.\\
(ii) If $f\ev0$, then
\[\ncP(v,T)\gsl N.\]
If all solutions are nondegenerate, then $\ncP(v,T)\gsl 2^{N-1}$.\eT

Similarly to \cite{FW89, CLZ90}, the aim of this paper is to establish additional rotational solutions for special
arrangement of masses $m_j$ and lengths $\ell_j$ instead of the nondegenerate assumption.
In the following result, we consider $v=(v_1,v_2,\cd,v_{\sss N})\in\Z^N_1$ with zero components. Let
\be N_0=\,^\#\big\{i~\big|~1\lsl i\lsl N,~v_i=0\big\}.\lb{eq:thm2.N0}\ee

\bT\lb{thm:Lsys.rot.new}\it Assume that $f$ satisfies \eqref{eq:intf=0} and
\be |f(t)|\lsl M_0,~\fa t\in\R.\ee
If $1\lsl N_0\lsl N-1$, there exist masses $m_j$, lengths $\ell_j$ ($1\lsl j\lsl N$), and constants
$T_1,T_2>0$ such that for each $T\in[T_1,T_2]$, we have\\
(i) If $f\nev0$, then
\[\ncP(v,T)\gsl (N-N_0+1)2^{N_0}.\]
(ii) If $f\ev0$, then
\[\ncP(v,T)\gsl (N-N_0)2^{N_0}.\]\eT




\sct{Variational settings}

For $x=(x_1,x_2,\cd,x_{\sss N})\in\cE$, denote by $||x||$ the standard $L^2$-norm of $x(t)$, i.e.,
\[||x||^2=\int^T_0|x(t)|^2\rd t=\sum^N_{i=1}\int^T_0x^2_i(t)\rd t.\]
The functional $\cL:\cE\ra\R$ defined by \eqref{eq:defcA} can be written as
\be\cL(x)=\cL_1(x)+\cL_2(x)+\cL_3(x),\lb{eq:cA=cA1+2+3}\ee
where
\be \cL_1(x)=\frac{1}{2}\<A(q)\dot q,\dot q\>,~\cL_2(x)=\int^T_0V(q)\rd t,~\cL_3(x)=\int^T_0f(t)\cdot q\rd t.\lb{eq:defcA123}\ee
According to \eqref{eq:cE=TNxE'}, we write $x\in\cE$ as
\[x=\bar x+\wtd x,\quad \bar x=(\bar x_1,\bar x_2,\cd,\bar x_{\sss N})\in T^N,
\quad \wtd x=(\wtd x_1,\wtd x_2,\cd,\wtd x_{\sss N})\in\wtd E.\]
Since the matrix $A(q)$ defined by \eqref{eq:defA(q)pp} is positive definite and periodic in $q$,
there exists a constant $\lm>0$ depending only on masses and lengths such that
\be\lm=\inf\lt\{\frac{A(q)p\cdot p}{|p|^2}~\big|~q\in T^N,\,p\in\R^N\bs\{0\}\rt\}.\lb{eq:eign.min}\ee
Let
\be \ga_1=2\pi^2\lt(\sum^N_{i=1}\aa_i\ell^2_iv^2_i+\sum_{i<j,\,v_i=v_j}2\aa_j\ell_i\ell_jv_iv_j\rt).\lb{eq:ga1}\ee 
We have

\bP\lb{prop:QVf}
\begin{align}
&\cL_1(x)\gsl\frac{\lm}{2}||\dot x||^2+T^{-1}2\pi^2|v|^2\lm,\lb{eq:Q>}\\
&\cL_1(\bar x)\lsl T^{-1}\ga_1,~\cL_2(\bar x)=T\sum_{i\in I_0}\bb_i\cos\bar x_i,
~\cL_3(\bar x)=\int^T_0f(t)\cdot\frac{2\pi vt}{T}\rd t=f_v,\lb{eq:cA123bx}\\
&\lt|\cL_2(x)-\cL_2(\bar x)\rt|\lsl\frac{T^{\frac{3}{2}}}{2\pi}\sum^N_{i=1}\bb_i||\dot x_i||,\lb{eq:Vq}\\
&\lt|\cL_3(x)-\cL_3(\bar x)\rt|\lsl \frac{T}{2\pi}\sum^N_{i=1}||f||\cdot||\dot x_i||
\lsl\frac{T^{\frac{3}{2}}}{2\pi}\sum^N_{i=1}M_0||\dot x_i||.\lb{eq:fq}\end{align}
\eP

\pf For \eqref{eq:Q>}, by \eqref{eq:defcA123} and \eqref{eq:eign.min}, we have
\[\cL_1(x)=\frac{1}{2}\<A(q)\dot q,\dot q\>_{L^2}\gsl\frac{\lm}{2}\lt|\lt|\dot x
+\frac{2\pi v}{T}\rt|\rt|^2=\frac{\lm}{2}||\dot x||^2+\frac{2\pi^2|v|^2\lm}{T}.\]

For \eqref{eq:cA123bx}, we claim that
\begin{align}
&\int^T_0\cos\lt(\bar x_i+2\pi v_i\frac{t}{T}\rt)\rd t
=\lt\{\bgd0,&&v_i\neq0,\\T\cos\bar x_i,&&v_i=0,\egd\rt.~~1\lsl i\lsl N,\lb{eq:cosv=0}\\
&\lt|\int^T_0\cos\lt(x_i+2\pi v_i\frac{t}{T}\rt)\rd t-\int^T_0\cos\lt(\bar x_i+2\pi v_i\frac{t}{T}\rt)\rd t\rt|
\lsl\frac{T\sqrt{T}}{2\pi}||\dot x_i||.\lb{eq:cosx-cosx^0}
\end{align}
In fact
\[\bgd &\lt|\int^T_0\cos\lt(x_i+2\pi v_i\frac{t}{T}\rt)\rd t-\int^T_0\cos\lt(\bar x_i+2\pi v_i\frac{t}{T}\rt)\rd t\rt|\\
&=\lt|\int^T_0\lt(\cos\lt(x_i+2\pi v_i\frac{t}{T}\rt)-\cos\lt(\bar x_i+2\pi v_i\frac{t}{T}\rt)\rt)\rd t\rt|\\
&=\lt|-2\int^T_0\sin\lt(\frac{\wtd x_i}{2}+\bar x_i+2\pi v_i\frac{t}{T}\rt)\sin\lt(\frac{\wtd x_i}{2}\rt)\rd t\rt|\\
&\lsl2\lt(\int^T_0\sin^2\lt(\frac{\wtd x_i}{2}+\bar x_i+2\pi v_i\frac{t}{T}\rt)\rd t\rt)^{\frac{1}{2}}
\lt(\int^T_0\sin^2\frac{\wtd x_i}{2}\rd t\rt)^{\frac{1}{2}}\\
&\lsl2\sqrt{T}\lt(\int^T_0\sin^2\frac{\wtd x_i}{2}\rd t\rt)^{\frac{1}{2}}
\lsl2\sqrt{T}\lt(\int^T_0\frac{\wtd x^2_i}{4}\rd t\rt)^{\frac{1}{2}}\\
&=\sqrt{T}||\wtd x_i||\lsl\frac{T\sqrt{T}}{2\pi}||\dot x_i||.\egd\]
The last inequality follows from the Wirtinger's inequality.

By \eqref{eq:defcA123}, \eqref{eq:defA(q)pp}, \eqref{eq:q=x+tv} and \eqref{eq:cosv=0}, we have
\[\bgd \cL_1(\bar x)&=\sum^N_{i=1}\aa_i\ell^2_i\cdot2\pi^2v^2_i\frac{1}{T}
+\sum_{1\lsl i<j\lsl N}\aa_j\ell_i\ell_j\cdot\frac{4\pi^2}{T^2}v_iv_j
\int^T_0\cos\lt(\bar x_i-\bar x_j+2\pi(v_i-v_j)\frac{t}{T}\rt)\rd t\\
&=\frac{2\pi^2}{T}\sum^N_{i=1}\aa_i\ell^2_iv^2_i+\sum_{i<j,\,v_i=v_j}
\aa_j\ell_i\ell_j\cdot\frac{4\pi^2}{T^2}v_iv_j\int^T_0\cos(\bar x_i-\bar x_j)\rd t\\
&\lsl\frac{2\pi^2}{T}\sum^N_{i=1}\aa_i\ell^2_iv^2_i+\frac{4\pi^2}{T}
\sum_{i<j,\,v_i=v_j}\aa_j\ell_i\ell_jv_iv_j=T^{-1}\ga_1.\egd\]
By \eqref{eq:defcA123}, \eqref{eq:q=x+tv}, \eqref{eq:defV} and \eqref{eq:intf=0}, we have
\[\bgd\cL_2(\bar x)&=\int^T_0V\lt(\bar x+2\pi v\frac{t}{T}\rt)\rd t\\
&=\int^T_0\sum^N_{i=1}\bb_i\cos\lt(\bar x_i+2\pi v_i\frac{t}{T}\rt)\rd t=T\sum_{i\in I_0}\bb_i\cos\bar x_i,\\
\cL_3(\bar x)&=\int^T_0f(t)\cdot\lt(\bar x+2\pi v\frac{t}{T}\rt)\rd t=\int^T_0f(t)\cdot\frac{2\pi vt}{T}\rd t.\egd\]

The inequality \eqref{eq:Vq} follows from \eqref{eq:defcA123}, \eqref{eq:defV} and \eqref{eq:cosx-cosx^0}.

For \eqref{eq:fq}, by \eqref{eq:defcA123}, \eqref{eq:q=x+tv} and \eqref{eq:intf=0}, we have
\[\bgd |\cL_3(x)-\cL_3(\bar x)|&=\lt|\int^T_0f\cdot q\rd t-\int^T_0f2\pi v\frac{t}{T}\rd t\rt|=\lt|\int^T_0f\cdot(\bar x+\wtd x)\rd t\rt|\\
&=\lt|\int^T_0f\cdot\wtd x\rd t\rt|=\lt|\int^T_0\sum^N_{i=1}f_i\cdot\wtd x_i\rd t\rt|\lsl\sum^N_{i=1}\int^T_0|f_i|\cdot|\wtd x_i|\rd t\\
&\lsl\sum^N_{i=1}||f_i||\cdot||\wtd x_i||\lsl\sum^N_{i=1}||f_i||\cdot\frac{T}{2\pi}||\dot x_i||\lsl\frac{T||f||}{2\pi}||\dot x||.\egd\]
The third equality follows from (\rf{eq:cA123bx}), and the last inequality follows from the Wirtinger's inequality.~\QED

By the above Proposition, we have
\be\lt|\cL_2(x)+\cL_3(x)-\big(\cL_2(\bar x)+\cL_3(\bar x)\big)\rt|
\lsl\frac{T^{\frac{3}{2}}}{2\pi}\sum^N_{i=1}(\bb_i+M_0)||\dot x_i||.\lb{eq:Vq+fq}\ee

\bP\lb{prop:cA.PS} \it For the functional $\cL$ defined by \eqref{eq:defcA}, we have\\
(i) $\cL$ is bounded from below.\\
(ii) $\cL$ satisfies the (PS) condition.\eP
\pf For (i), by , we have
\begin{align}\cL(x)&\gsl\frac{\lm}{2}||\dot x||^2+T^{-1}2\pi^2|v|^2\lm-T\sum^N_{i=1}\bb_i
-\frac{T^{\frac{3}{2}}||f||}{2\pi}||\dot x||+f_v\nn\\
&=\frac{\lm}{2}\sum^N_{i=1}\lt(||\dot x||-\frac{T^{\frac{3}{2}}||f||}{2\pi\lm}\rt)^2+a_0,~\fa x\in\cE,\lb{eq:cA.PS.1}\end{align}
where
\be a_0=T^{-1}2\pi^2|v|^2\lm+f_v-T\sum^N_{i=1}\bb_i-\frac{T^3||f||^2}{8\pi^2\lm}.\lb{eq:a0}\ee

For (ii), let $\{x_m\}$ be a sequence in $\cE$ such that
\be\cL(x_m)\lsl C~~\text{and}~~\cL'(x_m)\ra0~\text{as}~m\ra\ift.\lb{eq:LcA.ps}\ee
Then \eqref{eq:cA.PS.1} yields that
\be||\dot x_m||\lsl||x_m||_E\lsl C',\lb{eq:xnbdd}\ee
i.e., $\{x_m\}$ is bounded in $\cE$ since $T^N$ is compact. Without loss of generality,
assume that $x_m\rap x_0\in\cE$ and $x_m\ra x_0\in C(S^1,\R^N)$ as $m\ra\ift$, i.e.,
\be\lim_{m\ra\infty}||x_m-x_0||_{C^0}=0.\lb{eq:xm-x0.c0}\ee

Note that for arbitrary $x,y\in\cE$, by using the substitution \eqref{eq:q=x+tv} for simplicity, we have
\be\cL'(x)y=\cL'_1(x)y+\<V'(q)+f,y\>,\lb{eq:LcA'xy}\ee
where
\[\cL'_1(x)y=\<A(q)\dot q,\dot y\>+\lt\<\big(A'(q)y\big)\dot q,\dot q\rt\>,\]
\[\lt\<\big(A'(q)y\big)\dot q,\dot q\rt\>=\int^T_0\sum\aa_j\ell_i\ell_j\big(-\sin(q_i-q_j)(y_i-y_j)\big)\dot q_i\dot q_j\rd t.\]
Then
\[\bgd \cL'_1(x+y)y-\cL'_1(x)y&=\<A(q+y)(\dot q+\dot y),\dot y\>
+\lt\<(A'(q+y)y)(\dot q+\dot y),\dot q+\dot y\rt\>\\
&~~~~-\<A(q)\dot q,\dot y\>-\lt\<\big(A'(q)y\big)\dot q,\dot q\rt\>\\
&=\<A(q+y)\dot y,\dot y\>+\lt\<\big(A(q+y)-A(q)\big)\dot q,\dot y\rt\>\\
&~~~~+\lt\<\big(A'(q+y)y-A'(q)y\big)\dot q,\dot q\rt\>\\
&~~~~+2\lt\<\big(A'(q+y)y\big)\dot q,\dot y\rt\>
+\lt\<\big(A'(q+y)y\big)\dot y,\dot y\rt\>\\
&\gsl\lm||\dot y||^2+o\big(||y||\big).\egd\]
We have
\be \cL'_1(x_m)(x_m-x_0)-\cL'_1(x_0)(x_m-x_0)
\gsl\lm||\dot x_m-\dot x_0||^2+o(||x_m-x_0||).\lb{eq:Q'(x+y)y-Q'(x)y>}\ee
By \eqref{eq:LcA'xy} and \eqref{eq:Q'(x+y)y-Q'(x)y>}, we have
\[\bgd&\cL'(x_m)(x_m-x_0)-\cL'(x_0)(x_m-x_0)\\
&=\cL'_1(x_m)(x_m-x_0)+\<V'(q_m)+f,x_m-x_0\>\\
&~~~~-\cL'_1(x_0)(x_m-x_0)-\<V'(q_0)+f,x_m-x_0\>\\
&\gsl\lm||\dot x_m-\dot x_0||^2+o(||x_m-x_0||).\egd\]
Then \eqref{eq:LcA.ps} and \eqref{eq:xm-x0.c0} yield
$||\dot x_m-\dot x_0||\ra0$ as $m\ra\ift$. Hence $x_m\ra x_0$ in $\cE$. \QED\\

{\it Proof of \rT{thm:Lsys.rot}.} For (i), by \rP{prop:cA.PS} and Ljusternik-Schnirelman theory, we have
\[\ncP(v,T)\gsl \text{cat}(T^N)=N+1.\]
If all solutions are nondegenerate, by Morse inequalities, we have
\[\ncP(v,T)\gsl~\text{the sum of Betti numbers of}~T^N=2^N.\]

For (ii), by \rP{prop:cA.PS} and $S^1$-equivariant Ljusternik-Schnirelman theory, we have
\[\ncP(v,T)\gsl\text{cat}(T^N/S^1)=\text{cat}(T^{N-1})=N.\]
If all solutions are nondegenerate, by Morse inequalities, we have
\[\ncP(v,T)\gsl~\text{the sum of Betti numbers of}~T^{N-1}=2^{N-1}.~\QED\]

\sct{An abstract critical point result}

In \cite{CLZ90}, K.~C. Chang, Y.~Long and E.~Zehnder prove an abstract critical point result
(Theorem 3.1) and apply it to oscillatory solutions of \eqref{eq:Lsys}. In this section, we
establish corresponding result and apply it to rotational solutions. The difference is that
$N_0$ defined by \eqref{eq:thm2.N0} satisfies
\[N_0\lgd =N,&\quad \text{oscillatory case, i.e.}, v=0,\\ \lsl N-1,&\quad \text{rotational case, i.e.}, v\ne0.\rgd\]

The main idea is to find an $m$-dimensonal subtorus $T^m\sbs T^N$ with an $(m-1)$-dimensonal
subtorus $T^{m-1}\sbs T^m$ such that the inclusion
\[(T^{m},T^{m-1})~\hra~(I^b,I^a)\]
induces an injective map in homology and a surjective map in cohomology. By
Ljusternik-Schirelman theory, there exist at least $m$ critical points between the
two levels $a$ and $b$.

\begin{center}
\rotatebox{-0}{\includegraphics[width=8cm]{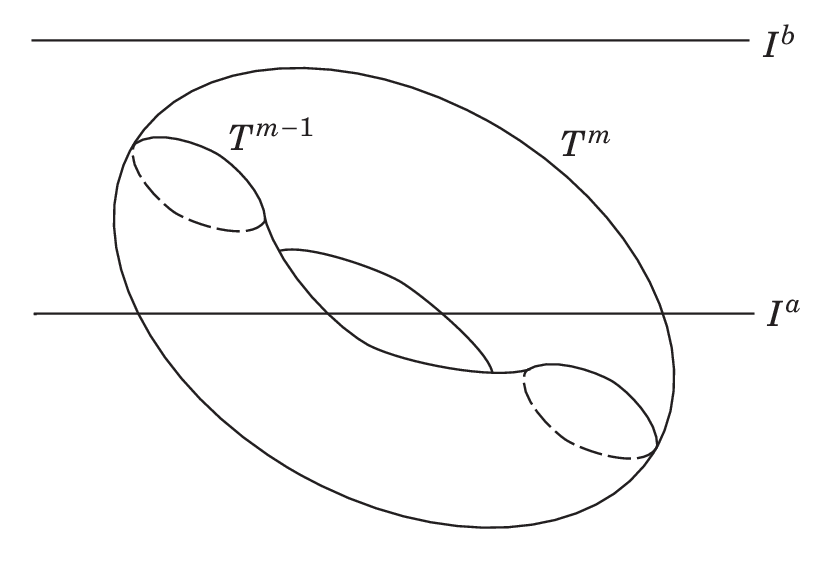}}
\end{center}

Let
\be n=2^{N_0}.\ee
For each $1\lsl k\lsl n$, represent $k-1$ by binary numbers as follows
\be k-1=\sum^{N_0}_{i=1}\tau_i(k)2^{ N_0-i},~\tau_i(k)\in\big\{0,1\big\},\lb{eq:bi.num}\ee
and set
\be z(k)=\lt(z_1(\tau_1(k)),z_2(\tau_2(k)),\cd,z_{\sss N_0}(\tau_{\sss N_0}(k))\rt)\in T^{N_0},\lb{eq:AbsCrit.zk}\ee
where $z_i:\big\{0,1\big\}\ra S^1,~1\lsl i\lsl N_0$, are injective maps.

Define an $(N-N_0)$-dimensional subtorus $S_k\sbs T^N$ by
\be S_k=\{z(k)\}\ts T^{N-N_0}.\lb{eq:Sk}\ee
For $k,j$ satisfying $1\lsl j\lsl k$ and
\[k-j=2^{N_0-m},~~\text{for some}~1\lsl m\lsl N_0,\]
define an $(N-N_0+1)$-dimensional subtorus $T_{k,\,j}\sbs T^N$ containing subtorus $S_k$
and $S_j$ as follows
\be T_{k,\,j}=\big\{(x_1,\cd,x_{\sss N_0})\in T^{N_0}~\big|~x_i=z_i(\tau_i(k)),
~i\neq m\big\}\ts T^{N-N_0}\cong T^{N-N_0+1}.\lb{eq:Tkk'}\ee

For $1\lsl k\lsl n,~1\lsl i\lsl N_0$, let
\begin{gather}
T^{N-1}_{k,\,i}=\big\{(y_1,\cd,y_{\sss N})~\big|~y_i=z_i(1-\tau_i(k))~\text{and}~
~y_j\in S^1~\text{for}~j\neq i\big\}\cong T^{N-1},\nn\\
Q_k=\bigcup^{N_0}_{i=1} T^{N-1}_{k,\,i},~M_k=\bigcap^{n}_{j=k}Q_j.\lb{eq:hQk,hMk}
\end{gather}

By the proof of Lemma 3.1 of \cite{CLZ90}, there is a strong deformation retraction
\be r_k: T^N\bs S_k\ra Q_k.\lb{eq:sdrk}\ee
Let
\begin{align}
W_k&=T^{N}\bs\lt[S_n\cup\lt(\bigcup^{n-1}_{i=k}r^{-1}_n\cd r^{-1}_{i+1}(S_i)\rt)\rt],\\
R_k&=r_k\circ r_{k+1}\circ\cd\circ r_n: \text{dom}(r_k\circ r_{k+1}\circ\cd\circ r_n)=W_k\,\ra T^{N}.\end{align}

\bP\lb{prop:CLZ90}\it Assume $I\in C^1(\cE,\R)$ satisfies (PS) condition and $\dsty a_0=\inf_{x\in\cE}I(x)>-\ift$.
If there exist real numbers $a_0<a_1<\cd<a_n$, where $n=2^{N_0}$ such that
\be I_{T^N\ts\{0\}}<a_n,~I|_{M_{k+1}\ts\{0\}}<a_k<I|_{p^{-1}R^{-1}_{k+2}(S_{k+1})},~1\lsl k\lsl n-1,\lb{eq:Abs.Crit}\ee
(Here we set $R^{-1}_{n+1}(S_n)=S_n$).
Denote by $\crit(I)$ the set of critical points of $I$, we have
\[\ncrit(I)\gsl (N-N_0+1)2^{N_0}.\]
If $I$ is $S^1$-equivariant, then
\[\ncrit(I)\gsl (N-N_0)2^{N_0}.\]\eP

\pf For brevity, let
\[m=N-N_0+1.\]
We claim that
\begin{align}
^\#\bigg(\crit(I)\cap \mrg I^{a_1}\bigg)&\gsl m,\quad k=1,\lb{eq:num.1}\\
^\#\bigg(\crit(I)\cap\big(\mrg I^{a_k}\bs\mrg I^{a_{k-1}}\big)\bigg)&\gsl m,\quad 2\lsl k\lsl n-1,\lb{eq:num.2}\\
^\#\bigg(\crit(I)\cap\big(\cE\bs I^{a_{n-1}}\big)\bigg)&\gsl m,\quad k=n,\lb{eq:num.3}
\end{align}
where
\[I^a=\{x\in\cE~|~I(x)\lsl a\}\text{\quad and\quad}\mrg{I}^a=\{x\in\cE~|~I(x)<a\}.\]
By \eqref{eq:num.1}-\eqref{eq:num.3}, we have
\[\ncrit(I)\gsl m+(n-2)m+m=mn=(N-N_0+1)2^{N_0}.\]

\nid{\bf Case 1.} $k=1$. By \eqref{eq:Abs.Crit}, there exists $\dl_1\in (0,a_1-a_0)$ such that
\[S_1\ts\{0\}\sbs I^{a_1-\dl_1}.\]
The Ljusternik-Schnirelman category of $I^{a_1-\dl_1}$ is
\[\text{cat}(I^{a_1-\dl_1})\gsl\text{cat}(S_1)=\text{cat}(T^{m-1})=m.\]
As $a_0=\inf I>-\ift$ and $I$ satisfies (PS), we conclude \eqref{eq:num.1}.

\nid{\bf Case 2.} $2\lsl k\lsl n-1$. We proceed in several steps.

{\bf Step 1.}  Choose $\dl\in(0,a_k-a_{k-1})$  sufficiently small, abbreviate
\[b=a_k-\dl,~a=a_{k-1} \]
and set
\[\bgd
X&=\cE\bs\lt[p^{-1}(S_n)\cup p^{-1}R^{-1}_n(S_{n-1})\cup\cd\cup p^{-1}R^{-1}_{k+2}(S_{k+1})\rt],\\
X'&=\cE\bs\lt[p^{-1}(S_n)\cup p^{-1}R^{-1}_n(S_{n-1})\cup\cd\cup p^{-1}R^{-1}_{k+1}(S_k)\rt],\\
Y&=I^b,~Y'=I^a,~Z=M_{k+1}\ts\{0\},~Z'=M_k\ts\{0\}.
\egd\]
By \eqref{eq:Abs.Crit}, we have
\[\ba{ccccc} X & \sps &   Y  & \sps &   Z\\
          \cup &      & \cup &      & \cup\\
            X' & \sps &   Y' & \sps &   Z'.\ea\]

{\bf Step 2.} Consider the following inclusion maps
\[(T_{k,\,k'},S_{k'}) \xra{\rho} (M_{k+1},M_k) \xra[x\mt (x,0)]{\bb} (Z,Z') \xra{j} (Y,Y').\]
By the proof on page 195 of \cite{CLZ90}, the continuous map
\[\psi=j\circ\bb\circ\rho: (T_{k,\,k'},S_{k'})\,\ra(I^b,I^a)\]
induces a monomorphism
\[\psi_\ast: H_\ast(T_{k,\,k'},S_{k'}) \ra H_\ast(I^b,I^a),\]
and an epimorphism
\[\psi^\ast: H^\ast(I^b) \ra H^\ast(T_{k,\,k'}).\]

{\bf Step 3.} We claim that
\[H_q(T_{k,\,k'},S_{k'})\cong  H_q\lt(T^m,T^{m-1}\rt)\neq0,~1\lsl q\lsl m.\]
In fact, the first isomorphism follows from \eqref{eq:Sk} and \eqref{eq:Tkk'}. Consider the
following long exact sequence
\[\cd \ra H_q(T^{m-1}) \mra{~i_q~} H_q(T^m ) \mra{~j_q~} H_q(T^m,T^{m-1}) \mra{~\pt_q~} H_{q-1}(T^{m-1}) \mra{~i_{q-1}~} \cd\]
We have
\[\bgd H_q(T^{m-1}) &\cong \text{Ker}i_q\op\text{Im}i_q \cong \Z^{C^q_{m-1}},\\
           H_q(T^m) &\cong \text{Ker}j_q\op\text{Im}j_q \cong \Z^{C^q_m},\\
   H_q(T^m,T^{m-1}) &\cong \text{Ker}\pt_q\op\text{Im}\pt_q,\\
   \text{Im}i_q&=\text{Ker}j_q,~\text{Im}j_q=\text{Ker}\pt_q,~\text{Im}\pt_q=\text{Ker}i_{q-1}.\egd\]
Hence, for $1\lsl q\lsl m$, the $q$-th Betti-number of the pair $(T^m,T^{m-1})$ satisfies
\[\bgd b_q(T^m,T^{m-1})&=\text{dim}\,H_q(T^m,T^{m-1})\\
&=C^q_m-C^q_{m-1}+\text{dim}\,\text{Ker}i_q+\text{dim}\,\text{Ker}i_{q-1}\\
&\gsl C^q_m-C^q_{m-1}\gsl1.\egd\]

{\bf Step 4.} Choose homology classes
\[0\neq[u_q]\in H_q(T_{k,\,k'},S_{k'}),~1\lsl q\lsl m,\]
and a cohomology class $0\neq\tht\in H^1(T_{k,\,k'})\cong H^1(T^m)$ such that
\[[u_{q+1}]\cap\tht=[u_q],~1\lsl q\lsl m-1.\]
Let
\[[\sg_q]=\psi_\ast[u_q],~1\lsl q\lsl m,\]
and choose
\[\om\in(\psi^\ast)^{-1}(\tht),\]
then
\[0\neq[\sg_q]\in H_q(I^b,I^a)~~\text{and}~~0\neq\om\in H^1(I^b).\]
By Proposition 5.6.16 (on page 254) of \cite{eSp66}, we have
\begin{align}
[\sg_{q+1}]\cap\om&=\big(\psi_\ast([u_{q+1}])\big)\cap\om=\psi_\ast([u_{q+1}\cap\psi^\ast\om])\nn\\
&=\psi_\ast([u_{q+1}\cap\tht])=\psi_\ast([u_q])=[\sg_q],~1\lsl q\lsl m-1.\nn\end{align}
By definition 1.1 (on page 10) of \cite{kcC93}, we have
\[[\sg_1]<[\sg_2]<\cd<[\sg_m].\]
Apply the Ljusternik-Schnirelman theory (see Corollary 3.3 (on page 106) of \cite{kcC93}),
we conclude \eqref{eq:num.2}.

\nid{\bf Case 3.} $k=n$. The proof is similar to that of Case 1. In Step 1, we set
\[X=Y=\cE,~Z=T^N\ts\big\{0\big\},~Z'=M_n\ts\{0\}=Q_n\ts\big\{0\big\},\]
and conclude $N-N_0+1$ critical points above the level $a_{n-1}$.

If $I$ is $S^1$-equivariant, then replace $X,Y,Z,X',Y',Z'$ by their quotient spaces under
the $S^1$-action. We have
\[\bgd \text{cat}(I^{a_1-\dl_1}/S^1)&\gsl\text{cat}(S_1/S^1)=\text{cat}(T^{m-2})=m-1,\\
H_q(T_{k,\,k'}/S^1,S_{k'}/S^1)&\cong H_q\lt(T^m/S^1,T^{m-1}/S^1\rt)\\ &\cong H_q\lt(T^{m-1},T^{m-2}\rt)\neq0,~1\lsl q\lsl m-1.\egd\]
Then repeat Step 4 in Case 2 by replacing $m$ by $m-1$.~\QED

For each homeomorphism
\be \phi: T^N\ra T^N,\lb{eq:phi}\ee
set
\begin{gather*}
\wtd S_k=\phi^{-1}\lt(S_k\rt)\cong T^{N-N_0},
~\wtd T_{k,\,j}=\phi^{-1}\lt(T_{k,\,j}\rt)\cong T^{N-N_0+1},\\
\wtd T^{N-1}_{k,\,i}=\phi^{-1}\lt(T^{N-1}_{k,\,i}\rt)\cong T^{N-1},
~\wtd Q_k=\phi^{-1}\lt(Q_k\rt),~\wtd M_k=\phi^{-1}\lt(M_k\rt),\\
\wtd r_k=\phi^{-1}\circ r_k,~~\wtd W_k=\phi^{-1}\lt(W_k\rt),
\wtd R_k=\phi^{-1}\circ\lt(R_k\rt)=\wtd r_k\circ \wtd r_{k+1}\circ\cd\circ \wtd r_n: \wtd W_k\ra T^N.\end{gather*}
By \eqref{eq:sdrk}, the map
\[\wtd r_k: T^N\bs \wtd S_k\ra \wtd Q_k,\]
is a strong deformation retraction. By \eqref{eq:hQk,hMk}, we have
\[\wtd Q_k=\bigcup^{N_0}_{i=1}T^{N-1}_{k,\,i},~\wtd M_k=\bigcap^n_{j=k}\wtd Q_j.\]

Since the above homeomorphism $\phi$ induces isomorphisms in homology and cohomology,
by \rP{prop:CLZ90}, we have

\bcr\lb{prop:thm23}\it Assume $I\in C^1(\cE,\R)$ satisfies (PS) condition and
$\dsty a_0=\inf_{x\in\cE}I(x)>-\ift$. If there exist real numbers $a_0<a_1<\cd<a_n,$
where $n=2^{N_0}$ such that
\be I_{T^N\ts\{0\}}<a_n,~I|_{\wtd M_{k+1}\ts\{0\}}<a_k
<I|_{p^{-1}\wtd R^{-1}_{k+2}(\wtd S_{k+1})},~1\lsl k\lsl n-1,\lb{eq:Abs.Crit}\ee
(Here we set $R^{-1}_{n+1}(S_n)=S_n$).
Denote by $\crit(I)$ the set of critical points of $I$, we have
\[\ncrit(I)\gsl (N-N_0+1)2^{N_0}.\]
If $\phi$ and $I$ are $S^1$-equivariant, then
\[\ncrit(I)\gsl (N-N_0)2^{N_0}.\]\ecr

In \eqref{eq:AbsCrit.zk}, we set
\be z(k)=\big(\pi(1-\tau_1(k)),\pi(1-\tau_2(k)),\cd,\pi(1-\tau_{\sss N_0}(k))\big).\lb{eq:zk=Pi(1-tauk)}\ee

The following result will be utilized in the proof of \rT{thm:Lsys.rot.new}.

\bP\lb{prop:Ok}\it For $2\lsl k\lsl 2^{\sss N_0}-1$, let
\be O_k=\lt\{(y_1,y_2,\cd,y_N)\in T^N~\Bigg|~y_i=0,~\text{if}~\tau_i(k)=1\rt\}.\ee
Then
\be S_k\sbs R^{-1}_{k+1}(S_k)\sbs O_k.\lb{eq:RSkOk}\ee\eP

\pf For induction on $N_0$ and brevity, let $m=N_0$,
\[\bgd & z^{(m)}(k)=z(k), && S^{(m)}_k=S_k, && O^{(m)}_k=O_k,\\
& Q^{(m)}_k=Q_k, && M^{(m)}_k=M_k, && r_{m,\,k}=r_k, && R_{m,\,k+1}=R_{k+1}.\egd\]
The map $r_{m,\,k}=r_k: T^N\bs S^{(m)}_k\ra Q^{(m)}_k$ is defined by
\[\lt(y^{(m)},y'\rt)\mt\lt(\sg y^{(m)}+(1-\sg)z^{(m)}(k),y'\rt),\]
where $\sg\gsl 1$ is uniquely determined by $k$ and $y^{(m)}$ (c.f. Lemma 3.1 of \cite{CLZ90}). We have
\begin{gather}
r^{-1}_{m,\,k}\lt(w^{(m)},\om_{m+1}\rt)=\lt\{\lm w^{(m)}+(1-\lm)z^{(m)}(k)~\bigg|~0<\lm\lsl 1\rt\}\ts\{\om_{m+1}\},\\
r^{-1}_{m,\,s+1}r^{-1}_{m,\,s}\lt(O^{(m)}_k\rt)
=r^{-1}_{m,\,s+1}\big(r^{-1}_{m,\,s}\lt(O^{(m)}_k\rt)\cap Q^{(m)}_{s+1}\big),
~~k+1\lsl s\lsl 2^m,\\
r^{-1}_{m,\,q}\cd r^{-1}_{m,\,p}(\om)\sbs
\lt\{\lm_q\cd\lm_p\om+\sum^{q-1}_{i=p}\lm_q\cd\lm_{i+1}(1-\lm_i)z^{(m)}(i)\rt.\nn\\
\lt.\hs{3cm}+(1-\lm_q)z^{(m)}(q)\Bigg|~0<\lm_i\lsl 1\rt\}\ts T^{N-m}.\lb{eq:rmarmb}\end{gather}

Note that $R_{m,\,k+1}\lt(S^{(m)}_k\rt)=S^{(m)}_k$. We claim that
\be R^{-1}_{m,\,s}\lt(O^{(m)}_k\rt)\sbs O^{(m)}_k,~~k\lsl 2^{m}-1,~~s=1,~k+1.\lb{eq:Rms}\ee

\begin{table*}[htbp]
\centering
\begin{tabular}{|c|c|c|}
\hline \multicolumn{1}{|c|}{} &
\multicolumn{1}{|c|}{$1\lsl k\lsl 2^{m-1}$} &
\multicolumn{1}{|c|}{$2^{m-1}+1\lsl k\lsl 2^m$}\\
\hline $z^{(m)}(k)$ & $\lt(\pi,z^{(m-1)}(k)\rt)$ & $\lt(0,z^{(m-1)}\lt(k-2^{m-1}\rt)\rt)$\\[1mm]
\hline $S^{(m)}_k$ & $\{\pi\}\ts S^{(m-1)}_k$ & $\{0\}\ts S^{(m-1)}_{k-2^{m-1}}$\\[1mm]
\hline $O^{(m)}_k$ & $S^1\ts O^{(m-1)}_k$ & $\{0\}\ts O^{(m-1)}_{k-2^{m-1}}$\\[1mm]
\hline $Q^{(m)}_k$ & $\lt(\{0\}\ts T^{N-1}\rt)\bigcup\lt(S^1\ts Q^{(m-1)}_k\rt)$
& $\lt(\{\pi\}\ts T^{N-1}\rt)\bigcup\lt(S^1\ts Q^{(m-1)}_{k-2^{m-1}}\rt)$\\[1mm]
\hline
\end{tabular}
\end{table*}

If $m=2$, we have $2\lsl k\lsl 2^2-1=3$,
\[\bgd &O^{(2)}_2=S^1\ts\{0\}\ts T^{N-2},\\
&O^{(2)}_3=\{0\}\ts S^1\ts T^{N-2}.\egd\]
Then
\[\bgd R^{-1}_{2,\,3}\lt(O^{(2)}_2\rt)&=r^{-1}_{2,\,4}r^{-1}_{2,\,3}\lt(O^{(2)}_2\rt)
=r^{-1}_{2,\,4}\Bigg(r^{-1}_{2,\,3}\lt(O^{(2)}_2\rt)\cap Q^{(2)}_4\Bigg)
\sbs r^{-1}_{2,\,4}\lt(O^{(2)}_2\rt)\sbs O^{(2)}_2,\\
R^{-1}_{2,\,1}\lt(O^{(2)}_2\rt)
&=r^{-1}_{2,\,4}r^{-1}_{2,\,3}r^{-1}_{2,\,2}r^{-1}_{2,\,1}\lt(O^{(2)}_2\rt)
=r^{-1}_{2,\,4}r^{-1}_{2,\,3}r^{-1}_{2,\,2}\Bigg(r^{-1}_{2,\,1}\lt(O^{(2)}_2\rt)\cap Q^{(2)}_2\Bigg)\\
&\sbs r^{-1}_{2,\,4}r^{-1}_{2,\,3}r^{-1}_{2,\,2}\lt(O^{(2)}_2\rt)
=r^{-1}_{2,\,4}r^{-1}_{2,\,3}\Bigg(r^{-1}_{2,\,2}\lt(O^{(2)}_2\rt)\cap Q^{(2)}_3\Bigg)\\
&\sbs r^{-1}_{2,\,4}r^{-1}_{2,\,3}\lt(O^{(2)}_2\rt)\sbs O^{(2)}_2,\\
R^{-1}_{2,\,4}\lt(O^{(2)}_3\rt)&=r^{-1}_{2,4}\lt(O^{(2)}_3\rt)\sbs O^{(2)}_3,\quad
R^{-1}_{2,\,1}\lt(O^{(2)}_3\rt)=r^{-1}_{2,\,4}r^{-1}_{2,\,3}r^{-1}_{2,\,2}r^{-1}_{2,\,1}\lt(O^{(2)}_3\rt)=\ept.\egd\]

Assume \eqref{eq:Rms} holds for $m=p$. We carry out the proof for $m=p+1$ as follows.

{\bf Case 1.} $2\lsl k\lsl 2^{m-1}-1$. We have $O^{(m)}_k=S^1\ts O^{(m-1)}_k$. Then
\[\bgd
r^{-1}_{m,\,2^{m-1}}\cd r^{-1}_{m,\,s}\lt(O^{(m)}_k\rt)
&\sbs S^1\ts R^{-1}_{m-1,\,s}\lt(O^{(m-1)}_k\rt)\sbs S^1\ts O^{(m-1)}_k,~~s=1,~k+1,\\
R^{-1}_{m,\,s}\lt(O^{(m)}_k\rt)
&=R^{-1}_{m,\,2^{m-1}+1}\lt(r^{-1}_{m,\,2^{m-1}}\cd r^{-1}_{m,\,s}\lt(O^{(m)}_k\rt)\rt)\\
&\sbs R^{-1}_{m,\,2^{m-1}+1}\lt(S^1\ts O^{(m-1)}_k\rt)\\
&\sbs S^1\ts\lt(O^{(m-1)}_k\cup R^{-1}_{m-1,\,1}\lt(O^{(m-1)}_k\rt)\rt)\sbs O^{(m)}_k.\egd\]

{\bf Case 2.} $k=2^{m-1}$. We have $O^{(m)}_{2^{m-1}}=S^1\ts O^{(m-1)}_{2^{m-1}}$. Then
\[\bgd R^{-1}_{m,\,2^{m-1}+1}\lt(O^{(m)}_{2^{m-1}}\rt)
&\sbs S^1\ts \lt(O^{(m-1)}_{2^{m-1}}\cup R^{-1}_{m-1,\,1}\lt(O^{(m-1)}_{2^{m-1}}\rt)\rt)\sbs O^{(m)}_{2^{m-1}},\\
R^{-1}_{m,\,1}\lt(O^{(m)}_{2^{m-1}}\rt)
&=R^{-1}_{m,\,2^{m-1}+1}\lt(r^{-1}_{m,\,2^{m-1}}\cd r^{-1}_{m,\,1}\lt(O^{(m)}_{2^{m-1}}\rt)\rt)\\
&\sbs R^{-1}_{m,\,2^{m-1}+1}\lt(S^1\ts R^{-1}_{m-1,\,1}\lt(O^{(m-1)}_{2^{m-1}}\rt)\rt)\\
&\sbs R^{-1}_{m,\,2^{m-1}+1}\lt(S^1\ts O^{(m-1)}_{2^{m-1}}\rt)\\
&\sbs S^1\ts\lt(O^{(m-1)}_{2^{m-1}}\cup R^{-1}_{m-1,\,1}\lt(O^{(m-1)}_{2^{m-1}}\rt)\rt)\sbs O^{(m)}_{2^{m-1}}.\egd\]

{\bf Case 3.} $2^{m-1}+1\lsl k\lsl 2^m-1$. We have $1\lsl k-2^{m-1}\lsl 2^{m-1}-1$,
$O^{(m)}_k=\{0\}\ts O^{(m-1)}_{k-2^{m-1}}$ and
\[r^{-1}_{m,\,s}\lt(O^{(m)}_k\rt)=\{0\}\ts r^{-1}_{m-1,\,s-2^{m-1}}\lt(O^{(m-1)}_{k-2^{m-1}}\rt),
~~k+1\lsl s\lsl 2^m.\]
Then
\[\bgd R^{-1}_{m,\,k+1}\lt(O^{(m)}_k\rt)
&=r^{-1}_{m,\,2^m}\cd r^{-1}_{m,\,k+1}\lt(O^{(m)}_k\rt)\\
&=\{0\}\ts R^{-1}_{m-1,\,k-2^{m-1}+1}\lt(O^{(m-1)}_{k-2^{m-1}}\rt)\\
&\sbs \{0\}\ts O^{(m-1)}_{k-2^{m-1}}=O^{(m)}_k,\\
R^{-1}_{m,\,1}\lt(O^{(m)}_k\rt)
&=R^{-1}_{m,\,2^{m-1}+1}\lt(r^{-1}_{m,\,2^{m-1}}\cd r^{-1}_{m,\,1}\lt(O^{(m)}_k\rt)\rt)\\
&=R^{-1}_{m,\,2^{m-1}+1}\lt(\{0\}\ts R^{-1}_{m-1,\,1}\lt(O^{(m-1)}_{k-2^{m-1}}\rt)\rt)\\
&\sbs R^{-1}_{m,\,2^{m-1}+1}\lt(\{0\}\ts O^{(m-1)}_{k-2^{m-1}}\rt)\\
&=\{0\}\ts R^{-1}_{m-1,\,1}\lt(O^{(m-1)}_{k-2^{m-1}}\rt)
\sbs\{0\}\ts O^{(m-1)}_{k-2^{m-1}}\sbs O^{(m)}_k.~\QED\egd\]

\sct{Proof of \rT{thm:Lsys.rot.new}}

\ssct{Restrictions of masses, lengths and periods}

Choose a permutation $\sg$ of $\big\{1,2,\cd,N\big\}$ such that
\[\bgt I_0=\big\{\sg(1),~\sg(2),~\cd,~\sg(N_0)\big\},\\
\sg(1)<\sg(2)<\cd<\sg(N_0).\egt\]
Define a homeomorphism $\phi: T^N\ra T^N$ by
\[(x_1,x_2,\cd,x_i,\cd,x_{\sss N})
\mt (x_{\sg(1)},x_{\sg(2)},\cd,x_{\sg(i)},\cd,x_{\sg(N)}).\]
Let
\begin{align}
\bb(k)&=\sum^{N_0}_{i=1}(-1)^{1-\tau_i(k)}\bb_{\sg(i)},~1\lsl k\lsl n=2^{\sss N_0},\lb{eq:thm2.bb(k)}\\
\ga&=\min_{1\lsl k\lsl n-1}(\bb(k+1)-\bb(k)),\lb{eq:thm2.ga}\\
\ga_2&=\frac{1}{8\pi^2\lm}\sum^N_{i=1}(\bb_i+M_0)^2,\lb{eq:thm2.ga2}\\
T_1&=\sqrt{\frac{\ga_1}{\ga}},~~T_2=\sqrt{\frac{\ga}{\ga_2}}.\lb{eq:thm2.T1T2}\end{align}

\bP\lb{prop:thm2.mlT} Assume that
\be\ga>\sqrt{\ga_1\ga_2}.\lb{eq:thm2.ga>}\ee
Then
\be\ga T^2>\ga_1+\ga_2T^4,~\fa T\in[T_1,T_2].\lb{eq:thm2.gaT2>}\ee\eP

\pf By \eqref{eq:thm2.T1T2} and \eqref{eq:thm2.ga>}, we have $T_1<T_2$. For each $T\in[T_1,T_2]$, we have
\[\ga T^2>\max\lt\{\ga_1,~\ga_2T^4\rt\}.~\QED\]

\bP\lb{prop:thm2.bb(k)GaOm} \it Let
\begin{align}
\Ga&=\lt\{\frac{1}{2}\big(\bb(k+1)-\bb(k)\big)~\Bigg|~1\lsl k\lsl n-1\rt\}.\lb{eq:thm2.Ga}
\end{align}
Then
\begin{align*}
\Ga&=\big\{\bb_{\sg(N_0)}\big\}\bigcup\lt\{\bb_{\sg(i)}-\sum^{N_0}_{j=i+1}\bb_{\sg(j)}~\Bigg|~1\lsl i\lsl N_0-1\rt\}.
\end{align*}\eP

\pf Let
\[u(k)=(1-\tau_1(k),\cd,1-\tau_{\sss N_0}(k)).\]

{\bf Case 1.} $N_0=1$. $n=2^{N_0}=2$. We have $k=1$,
\[\bgd &u(2)=0, &&\bb(2)=\bb_1,\\ &u(1)=1, &&\bb(1)=-\bb_1.\egd\]
Then $\Ga=\big\{\bb_1\big\}$.

{\bf Case 2.} $N_0\gsl2,~n=2^{N_0}\gsl4$.
\[\bgd u(k+1)&=(u',0,\udb{1,\cd,1}_{i}),\\ u(k)&=(u',1,\udb{0,\cd,0}_{i}),\egd\quad 0\lsl i\lsl N_0-1.\]
Then
\[\frac{1}{2}\big(\bb(k+1)-\bb(k)\big)=\lt\{\bgd \bb_{\sg(N_0)},&\quad i=0,\\ \bb_{\sg(N_0-i)}-\sum^{N_0}_{j=N_0-i+1}\bb_{\sg(j)},&\quad i=1,\cd,N_0-1.~\QED\egd\rt.\]


\ssct{Some estimates on potential energy part}

\bP\lb{prop:V.Mk<V.Sk}\it For $1\lsl k\lsl 2^{\sss N_0}-1$, we have
\be V|_{\wtd M_{k+1}}\lsl\bb(k)<\bb(k+1)\lsl V|_{\wtd R^{-1}_{k+2}(\wtd S_{k+1})}.\lb{eq:VMk<VRSk}\ee\eP

\pf For $1\lsl s\lsl m$, let
\be V_s\lt(\xi_{m-s+1}\cd,\xi_m\rt)=\sum^m_{j=m-s+1}\bb_{\sg(j)}\cos\xi_j,\ee
and
\be h_s(k)=\sum^m_{j=m-s+1}(-1)^{1-\tau_j(k)}\bb_{\sg(j)}.\ee
Then
\[V_m|_{M^{(m)}_k}=V|_{\wtd M_k},\quad
V_m\big|_{R^{-1}_{m,\,k+1}\lt(S^{(m)}_k\rt)}
=V\big|_{{\hat R^{-1}_{k+1}\lt(\wtd S_k\rt)}},\quad h_m(k)=\bb(k)\]
and
\be h_m(k)-h_{m-1}(k)=(-1)^{1-\tau_1(k)}\bb_{\sg(1)}.\ee

By \eqref{eq:RSkOk}, we have
\[V_m\big|_{R^{-1}_{m,\,k+1}\lt(S^{(m)}_k\rt)}\gsl V_m\big|_{O^{(m)}_k}\gsl h_m(k).\]
We claim that
\be V_m|_{M^{(m)}_k}\lsl h_m(k-1).\ee


\begin{table*}[htbp]
\centering
\begin{tabular}{|c|c|c|}
\hline \multicolumn{1}{|c|}{} &
\multicolumn{1}{|c|}{$M^{(m)}_k$} &
\multicolumn{1}{|c|}{$V_m\big|_{M^{(m)}_k}\lsl$}\\[1mm]
\hline $2\lsl k\lsl 2^{m-1}$
       & $\{\pi\}\ts M^{(m-1)}_k$
       & $-\bb_{\sg(1)}+V_{m-1}|_{M^{(m-1)}_k}$\\[1mm]
\hline $k=2^{m-1}+1$
       & $\{\pi\}\ts T^{N-1}$
       & $\dsty-\bb_{\sg(1)}+\sum^{m}_{j=2}\bb_{\sg(j)}=\bb\lt(2^{m-1}\rt)$\\[1mm]
\hline $2^{m-1}+2\lsl k\lsl 2^m$
       & $\lt(\{\pi\}\ts T^{N-1}\rt)\bigcup\lt(S^1\ts M^{(m-1)}_{k-2^{m-1}}\rt)$
       & $\bb_{\sg(1)}+V_{m-1}|_{M^{(m-1)}_{k-2^{m-1}}}$\\[1mm]
\hline
\end{tabular}
\end{table*}

{\bf Case 1.} $2\lsl k\lsl 2^{m-1}$. Assume that $V_{m-1}|_{M^{(m-1)}_k}\lsl h_{m-1}(k-1)$. Then
\[V_m|_{M^{(m)}_k}\lsl-\bb_{\sg(1)}+V_{m-1}|_{M^{(m-1)}_k}\lsl -\bb_{\sg(1)}+h_{m-1}(k-1)=h_m(k-1).\]

{\bf Case 2.} $k=2^{m-1}+1$. We have
\[\bgt (k-1)-1=2^{m-1}-1=0\cdot 2^{m-1}+1\cdot2^{m-2}+\cd+1\cdot2^0,\\
\big(\tau_1(k-1),\tau_2(k-1),\cd,\tau_m(k-1)\big)=(0,1,1,\cd,1),\\
h_{m-1}(k-1)=\sum^m_{j=2}(-1)^{1-\tau_j(k-1)}\bb_{\sg(j)}=\sum^m_{j=2}\bb_{\sg(j)}.\egt\]
Then $\dsty V_m|_{M^{(m)}_k}\lsl-\bb_{\sg(1)}+\sum^m_{j=2}\bb_{\sg(j)}=-\bb_{\sg(1)}+h_{m-1}(k-1)=h_m(k-1)$.

{\bf Case 3.} $2+2^{m-1}\lsl k\lsl 2^m$. We have
\[(k-1)-1=2^{m-1}+\lt[\lt(k-1-2^{m-1}\rt)-1\rt],~~0\lsl\lt(k-1-2^{m-1}\rt)-1\lsl 2^{m-1}-2.\]
Then
\[\tau_j(k-1)=\lt\{\bgd 1,&\quad j=1,\\ \tau_j\lt(k-1-2^{m-1}\rt),&\quad 2\lsl j\lsl m.\egd\rt.\]
Assume that
\[V_{m-1}|_{M^{(m-1)}_{k-2^{m-1}}}\lsl h_{m-1}\lt(k-2^{m-1}-1\rt).\]
We have
\[\bgd
V_m\big|_{S^1\ts M^{(m-1)}_{k-2^{m-1}}}&\lsl \bb_{\sg(1)}+V_{m-1}|_{M^{(m-1)}_{k-2^{m-1}}}\\
&\lsl \bb_{\sg(1)}+h_{m-1}\lt(k-1-2^{m-1}\rt)\\
&=\bb_{\sg(1)}+h_{m-1}(k-1)=h_m(k-1).\egd\]
Note that
\[\bgd V_m\big|_{\{\pi\}\ts T^{m-1}}
&\lsl -\bb_{\sg(1)}+\sum^m_{j=2}\bb_{\sg(j)}\lsl \bb_{\sg(1)}-\sum^m_{j=2}\bb_{\sg(j)}\\
&\lsl \bb_{\sg(1)}+h_{m-1}(k-1)=h_m(k-1).\egd\]
Hence $V_m|_{M^{(m)}_k}\lsl h_m(k-1)$.~\QED


\ssct{Proof of \rT{thm:Lsys.rot.new}}

For $\ga_1$ define by \eqref{eq:ga1} and $\bb(k)$ defined by \eqref{eq:thm2.bb(k)}, let
\begin{align}
C_1(k)&=T^{-1}\ga_1+T\bb(k)+f_v,~~1\lsl k\lsl n-1,\lb{eq:thm2.C1k}\\
C_2(k)&=T^{-1}2\pi^2|v|^2\lm+T\bb(k+1)+f_v-T^3\ga_2,\lb{eq:thm2.C2k}\end{align}
where
\[f_v=\int^T_0f(t)\cdot\frac{2\pi vt}{T}\rd t.\]
By \eqref{eq:thm2.gaT2>}, we have
\begin{align}C_2(k)-C_1(k)&>T(\bb(k+1)-\bb(k))-T^{-1}\ga_1-T^3\ga_2\nn\\
&\gsl T\cdot\ga-T^{-1}\ga_1-T^3\ga_2>0.\lb{eq:thm2.C2k>C1k}\end{align}
Let
\begin{align} a_n&=T^{-1}\ga_1+T\sum^{N_0}_{i=1}\bb_i+f_v+1,\\
a_k&=\frac{1}{2}(C_1(k)+C_2(k)),~~~~1\lsl k\lsl n-1.\lb{eq:thm2.ak}\end{align}
By \eqref{eq:thm2.C2k>C1k}, we have
\be C_1(k)<a_k<C_2(k).\lb{eq:thm2.C1k<ak<C2k}\ee
We claim that
\be a_1<a_2<\cd<a_n.\lb{eq:a1<a2<an}\ee
and
\be\cL|_{T^N\ts\{0\}}<a_n,~\cL|_{\wtd M_{k+1}\ts\{0\}}<a_k<\cL|_{p^{-1}\wtd R^{-1}_{k+2}(\wtd S_{k+1})},~1\lsl k\lsl n-1,\lb{eq:cAMk<ak<cASk}\ee

{\bf Step 1.} In fact, by \eqref{eq:thm2.ak} and \eqref{eq:thm2.C2k>C1k}, we have
\[\bgd a_{k+1}-a_k&=\frac{1}{2}(C_2(k+1)-C_2(k)+C_1(k+1)-C_1(k))\\
&=\frac{T}{2}(\bb(k+2)-\bb(k+1)+\bb(k+1)-\bb(k))\\
&\gsl\frac{T}{2}(\ga+\ga)=T\ga>0.\egd\]

{\bf Step 2.} We prove $\cL|_{\wtd M_{k+1}\ts\{0\}}<a_k$. In fact, by the definition of
$\cL_2$, we have
\[\cL_2|_{\wtd M_{k+1}\ts\{0\}}=TV|_{\wtd M_{k+1}},~~\cL_2|_{p^{-1}\wtd R^{-1}_{k+2}(\wtd S_{k+1})}=TV|_{\wtd R^{-1}_{k+2}(\wtd S_{k+1})}.\]
By \rP{prop:QVf} and \rP{prop:Ok}, we have
\[\cL_1|_{\wtd M_{k+1}\ts\{0\}}\lsl T^{-1}\ga_1,~\cL_2|_{\wtd M_{k+1}\ts\{0\}}\lsl T\bb(k),~\cL_3|_{\wtd M_{k+1}\ts\{0\}}=f_v.\]
Then
\[\cL|_{\wtd M_{k+1}\ts\{0\}}\lsl T^{-1}\ga_1+T\bb(k)+f_v=C_1(k)<a_k.\]

{\bf Step 3.} We prove $\cL|_{p^{-1}\wtd R^{-1}_{k+2}(\wtd S_{k+1})}>a_k$. In fact, by \rP{prop:Ok}, we have
\[\cL_2|_{p^{-1}\wtd R^{-1}_{k+2}(\wtd S_{k+1})}\gsl T\bb(k+1).\]
For each $x\in p^{-1}\wtd R^{-1}_{k+2}(\wtd S_{k+1})$, by \rP{prop:QVf},
we have
\[\bgd\cL(x)&\gsl\sum^N_{i=1}\lt[\frac{\lm}{2}||\dot x_i||^2
-\frac{T\sqrt{T}}{2\pi}(\bb_i+M_0)||\dot x_i||\rt]
+T^{-1}2\pi^2|v|^2\lm+\cL_2(\bar x)+f_v\\
&\gsl-\sum^N_{i=1}\frac{\lt(\frac{T\sqrt{T}}{2\pi}(\bb_i+M_0)\rt)^2}{2\lm}
+T^{-1}2\pi^2|v|^2\lm+T\bb(k+1)+f_v\\
&=-T^3\ga_2+T^{-1}2\pi^2|v|^2\lm+T\bb(k+1)+f_v=C_2(k)>a_k.~\QED\egd\]



\begin{thebibliography}{10}

\bibitem{Arn78}
VI.~Arnold,
{Mathematical Methods of Clasical Mechanics},
New York, Springer-Verlag Inc. 1978. 75-93.

\bibitem{BBG06}
J.~Bellazzini, V.~Benci, M.~Ghimenti,
{Periodic orbits of a one-dimensional non-autonomous Hamiltonian system},
{\it J. Differential Equations}~{\bf 230} (2006) 275-294.

\bibitem{CFS87}
A.~Capozzi, D.~Fortunato, A.~Salvatore,
{Periodic solutions of Lagrangian systems with bounded potential},
{\it J Math Anal Appl}~{\bf 124} (1987) 482-494.

\bibitem{CLZ90}
K.~C. Chang, Y.~Long, E.~Zehnder,
{Forced oscillation for the triple pendulum},
{\it Analysis, et cetera}, Academic Press, Boston, 1990.

\bibitem{kcC93}
K.~C. Chang,
{\it Infinite Dimensional Morse Theory and Multiple Solution Problems},
Birkh\"auser, Boston, 1993.

\bibitem{pF92JDE}
P.~Felmer, {Periodic solutions of spatially periodic Hamiltonian systems},
{\it J. Diff. Eqns.}~{\bf 98} (1992) 143--168.

\bibitem{pF92TMA}
P.~Felmer,
{Rotation type solutions for spatially periodic Hamiltonian systems},
{\it Nonlinear Analysis, T.M.A.}~{\bf 19} (1992) 409--425.

\bibitem{FW89}
G.~Fournier, M.~Willem,
{Multiple solutions of the forced double pendulum equation},
{\it Ann. Inst. H. Poincar\'e Anal. Non Lin\'eaire}~{\bf 6} (1989) 259--281.

\bibitem{cG91}
C.~Gol\'e,
{Monotone maps of $T^n\ts\R^n$ and their periodic orbits},
The geometry of Hamiltonian systems (ed. by T.~Ratin),
MSRI Publications 22 (Springer, New York) (1991) 341--366.

\bibitem{fJ94PLMS}
F.~Josellis,
{Ljusternik-Schnirelman theory for flows and periodic orbits for Hamiltonian systems on $T^n\ts\R^n$},
{\it Proc. London Math. Soc.}~{\bf 68} (1994) 641--672.

\bibitem{fJ94JDE}
F.~Josellis,
{Morse theory for forced oscillations of Hamiltonian systems on $T^n\ts\R^n$},
{\it J. Diff. Eqns.}~{\bf 2} (1994) 360--384.

\bibitem{rPa66}
R.~Palais,
{The Lusternik-Schnirelman theory on Banach manifolds},
{\it Topology}~{\bf 5}(1966), 115-132.

\bibitem{hQ}
H. Qiao,
{Rotational solutions for Hamiltonian systems on $\R^n\ts T^n$}.

\bibitem{pR88}
P. Rabinowitz,
{On a class of functionals invariant under a  action,}
{\it Trans. Am. math. Soc.}~{\bf 310} (1988) 303--311.

\bibitem{RP01}
Paolo Roselli,
{A multiplicity result for the periodically forced N-pendulum with nonzero-mean valued forcings},
{\it Nonlinear Anal. T.M.A.}~{\bf 43} (2001) 1019-1041.

\bibitem{eSp66}
E.~Spanier, {\it Algebraic Topology},
McGraw-Hill, New York, 1966.

\bibitem{Ta90}
G.~Tarantello,
{Multiple forced oscillations for the N-pendulum equation},
{\it Commun. Math. Phys.}~{\bf 132} (1990) 499--517.

\end{thebibliography}
\end{document}